\documentclass[a4paper,11pt,leqno]{amsart}
\usepackage{etex}
\usepackage[a4paper,left=25mm,right=25mm,top=30mm,bottom=30mm,marginpar=25mm]{geometry}
\usepackage[utf8x]{inputenc}
\usepackage[english]{babel}
\usepackage{amsmath}
\usepackage{amssymb}
\usepackage{amsthm}
\usepackage{enumerate}
\usepackage{empheq}
\usepackage{fancybox}
\usepackage{xcolor}
\usepackage{mathrsfs}
\usepackage{slashed}
\usepackage{tikz}
\usepackage{harpoon}
\usepackage{setspace}
\usepackage{fancybox}
\usepackage{esint}
\usepackage{bm}
\usepackage{hyperref}
\usepackage{subfigure}

\newcommand{\res}{\mathop{\hbox{\vrule height 7pt width .5pt depth 0pt
\vrule height .5pt width 6pt depth 0pt}}\nolimits}

\newcommand{\R}{\mathbb{R}}

\newcommand{\N}{\mathbb{N}}
\newcommand{\M}{\mathcal{M}}
\newcommand{\Ha}{\mathcal{H}}
\newenvironment{dimo}{\textit{Proof.\ }}{\begin{flushright}$\square$\end{flushright}}
\def\mint
   {\mkern12mu\hbox{\vrule height4pt
   depth-3.2pt width5pt}%
   \mkern-16.5mu\int}

\newtheorem*{thmnonumber}{Theorem}
\newtheorem{thm}{Theorem}[section]
\newtheorem{defi}[thm]{Definition}

\newtheorem{prop}[thm]{Proposition}
\newtheorem{conj}[thm]{Conjecture} 
\newtheorem{cor}[thm]{Corollary}
\newtheorem{lemma}[thm]{Lemma}
\newtheorem{rmk}[thm]{Remark}
\newtheorem{ackn}{Acknowledgements\!}

\title[A Convex Decomposition Formula for the Mumford-Shah Functional]{A Convex Decomposition Formula for the Mumford-Shah Functional in Dimension One}
\author{Marcello Carioni}
\address{Max Planck Institute for Mathematics in the Science, Inselstrasse 22, 04103 Leipzig, Germany
\vspace*{1mm}\\
Universit\"at W\"urzburg, Emil-Fischer-Stra{\ss}e 40, 97074 W\"urzburg, Germany
\vspace*{1mm}\\}
\email{marcello.carioni@mathematik.uni-wuerzburg.de}
\begin{document}
\maketitle
\begin{abstract}
We study the convex lift of Mumford-Shah type functionals in the space of rectifiable currents and we prove a convex decomposition formula in dimension one, for finite linear combinations of SBV graphs. We use this result to prove the equivalence between the minimum problems for the Mumford-Shah functional and the lifted one and, as a consequence, we obtain a weak existence result for calibrations in one dimension.
\end{abstract}
\begin{center}
%\vspace*{1mm}
\end{center}
\begin{center}
%\vspace*{1mm}
\end{center}

\section{Introduction}
The Mumford-Shah functional is one of the most important variational model for image segmentation. It was introduced in the late 80's by Mumford and Shah (\cite{BDB},\cite{OAB}) and it can be defined in its general form as
\begin{equation}\label{mumford}
J(u,K) = \int_{\Omega\setminus K} |\nabla u|^2 \, dx +  \beta \mathcal{H}^{n-1}(K) + \alpha \int_{\Omega\setminus K} |u-g|^2 \, dx,
\end{equation} 
where $\Omega \in \mathbb{R}^n$ is open, $K\subset \Omega$ is closed and such that $\mathcal{H}^{n-1}(K)<\infty$, $g\in L^{\infty} (\Omega)$, $u\in W^{1,2}(\Omega\setminus K)$ and $\beta$ and $\alpha$ are tuning parameters. \\
The idea of the model is that given $g$ representing the level of gray of an image, it is possible to get a ``smoother'' version of it, ``close'' to the starting one in the $L^2$ norm, by finding a minimizer of \eqref{mumford}. The gain of smoothness for the minimizers comes from penalizing the oscillation of the competitors (i.e. the Dirichlet energy) and the length of the contour, in order to avoid fractal behaviour of the boundary of the processed image.\\
The existence of minimizers for \eqref{mumford} was proved in \cite{ETF2} introducing a weak formulation obtained considering $u\in SBV(\Omega)$ and replacing the set $K$ with $S_u$, i.e. the singular set of $u$:
\begin{equation}\label{mum}
F(u) = \int_{\Omega} |\nabla u|^2 \, dx +  \beta\mathcal{H}^{n-1}(S_u) + \alpha\int_{\Omega} |u-g|^2 \, dx.
\end{equation}
It is worth to remark that when $\alpha = 0$ and $\beta=1$, $F$ is called homogeneous Mumford-Shah functional.\\
In the following years there have been a huge effort in understanding the regularity properties of the functional defined above. We can cite some relevant papers in this direction like \cite{HIG}, \cite{OTE}, \cite{CFM}. However, despite all the effort, the main conjecture proposed by Mumford and Shah in their seminal paper still remains open in its full generality.
\begin{conj}[Mumford, Shah]
Let $(u,K)$ be a pair minimizing \eqref{mum}. Then $K$ is locally union of finitely many $C^{1,1}$ embedded arcs.
\end{conj}
As pointed out for the first time in \cite{OTE}, a blow up limit of appropriate sequences of minimizers of \eqref{mum} is a local minimizer of the homogeneous Mumford-Shah functional; for this reason the characterization of these minimizers is directly related to the solution of the conjecture stated above. For example it is known that harmonic functions are local minimizers of \eqref{mum} (for $\alpha=0$ and $\beta=1$) in small domains and that the same result holds for step functions and triple junctions (\cite{TCM}). Moreover the main achievement in this direction is contained in \cite{CIS} and it answers affirmatively to a conjecture proposed by De Giorgi in \cite{PCD}:
\begin{equation}\label{crack}
u(\rho,\theta) = \sqrt{\frac{2\rho}{\pi}} \sin\left(\frac{\theta}{2}\right) \quad \rho>0,\ -\pi<\theta<\pi
\end{equation}
is a global minimizer of the homogeneous Mumford-Shah functional. \eqref{crack} is usually called crack-tip.

In \cite{TCM} Alberti, Bouchitt\'e and Dal Maso introduced the notion of calibration for the Mumford-Shah functional that resembles closely the classical theory for minimal surfaces by Harvey and Lawson (\cite{CGE}). With this technique, in \cite{TCM}, they were able to prove the minimality of some candidates for the homogeneous Mumford-Shah functional like the triple junction or reproving the minimality of harmonic functions in a very elegant way. However it remains open the problem of finding a calibration for the crack-tip and for general minimums in higher dimensions. 
It is therefore a relevant issue to understand if, given $u$ a minimum for the Mumford-Shah functional, then there exists a calibration for $u$.\\
This is the question we are going to address in this paper. Existence of calibration is a common issue also in the field of minimal surfaces and also there it is not completely solved. One can refer to the work of Federer \cite{RFC} for the classical results in this theory.\\
As for the Mumford-Shah functional the main result in this direction was obtained by Chambolle in \cite{CRF}. He proved the existence of a calibration in dimension one in a weak asymptotic sense using the following representation formula introduced in \cite{TCM}:
\begin{displaymath}
F(u)= \sup_{\phi\in K} \int_{\Gamma_u} \langle \phi, \nu_{\Gamma_u}\rangle \, d \mathcal{H}^{n}= \sup_{\phi\in K} \int_{\Omega\times \mathbb{R}} \langle\phi, D\textbf{1}_{\{u>t\}}\rangle, 
\end{displaymath}
where $K$ is the set of Borel vector fields $\phi:\Omega \times \R \rightarrow \R^{n+1}$ such that
\begin{equation}
\left\{\begin{array}{ll}
\phi^t(x,t) \geq \dfrac{|\phi^x(x,t)|^2}{4} -\beta(t-g)^2 & \forall x,t\\
\left|\displaystyle\int_{t_1}^{t_2} \phi^x(x,s)\,ds\right|\leq\alpha\ &\forall x,t_1,t_2.
\end{array}
\right.
\end{equation}
More precisely this representation formula is the particular case of a general one for ``local'' functionals in BV presented by Bouchitt\'e in \cite{CAN}.\\
In particular one can lift $F$ to higher dimension to obtain a convex functional $\mathcal{F}$ defined as
\begin{displaymath}
\mathcal{F}(w) = \sup_{\phi\in K\cap C_0} \int_{\Omega\times \mathbb{R}} \langle\phi , Dw\rangle
\end{displaymath}
for $w \in SBV(\Omega \times \R)$ decreasing in the last variable. If one is able to prove that given $u$ a minimizer of $F$, then $\textbf{1}_{\{u>t\}}$ is a minimizer of $\mathcal{F}$, then this would imply the existence of a calibration in a weak asymptotic sense by argument of convex analysis. Moreover, another important consequence is that one can compute the minimum of $F$ using the functional $\mathcal{F}$ that, being convex, allows for an efficient gradient descent method (\cite{AAF}).
Chambolle, in \cite{CRF}, was able to prove these facts in dimension one and he pointed out that the same results could be obtained building up a coarea-type formula for the previous functional generalising the classical coarea formula for functionals (\cite{CLD}, \cite{GCF}): 
\begin{equation}\label{solast}
\mathcal{F}(w) = \int_0^1 \mathcal{F}(\textbf{1}_{\{w(x,t)>s\}})\, ds,
\end{equation}
that is false for $\mathcal{F}$ as the example below shows:
\begin{displaymath}
u_1(x)=\left\{
\begin{array}{ll}
0 & \mbox{ if } x\leq 1/2\\
x & \mbox{ if } x> 1/2,
\end{array}
\right.
\qquad
u_2(x)=\left\{
\begin{array}{ll}
x & \mbox{ if } x\leq 1/2\\
1 & \mbox{ if } x> 1/2
\end{array}
\right.
\end{displaymath}
and 
$w(x,t) = (1/2)\textbf{1}_{\{u_1(x)> t\}} + (1/2)\textbf{1}_{\{u_2(x)> t\}}$.

\vspace{3mm}

In this article we use an alternative representation of the Mumford-Shah functional by rectifiable currents of the type
\begin{displaymath}
G(T)=\sup_{\phi\in K} \int_{\mathcal{M}} \theta\langle \phi , \nu_T \rangle d\mathcal{H}^{n},  
\end{displaymath}
where $T=(\mathcal{M},\xi,\theta)$ is a rectifiable current and $\nu_T$ is the normal to $\mathcal{M}$, and we start to exploit the validity of a general coarea-type formula for the functional $G$. In Section \ref{mais} we study the structure of the functional and we prove the following convex decomposition formula for a finite linear combination of graphs in dimension one. 
\begin{thmnonumber}[Convex decomposition formula]
Let $I$ be an open interval. Given $T= \sum_{i=1}^k \lambda_i\Gamma_{u_i}$ with $u_i\in SBV(I)$ and $\lambda_i>0$ such that $|\bigcup S_{u_i}| < +\infty$ there exists $k' \in \mathbb{N}$, $\{\mu_i\}_{i=1\ldots k'} > 0$ and $\{w_i\}_{i=1\ldots k'} \subset SBV(I)$ such that $T = \sum_{i=1}^{k'} \mu_i\Gamma_{w_i}$ and 
\begin{displaymath}
G(T) =  \sum_{i=1}^{k'} \mu_i G(\Gamma_{w_i}).
\end{displaymath}
\end{thmnonumber}
The previous formula can be viewed as a variant of the generalized coarea formula in the sense of \eqref{solast}, when the latter is applied to finite linear combination of graphs. \\
The immediate consequence of this result is the following theorem that links the minimizers of \eqref{mum} with the minimizers of $G$:
\begin{thmnonumber}
Given $u\in SBV(I)$ a minimizer of the Mumford-Shah functional, $\Gamma_u$ (i.e. the graph associated to $u$) is a minimizer of $G$ among all the linear combinations of graphs of the form $T = \sum_{i=1}^k \lambda_i \Gamma_{u_i}$ with $\partial \Gamma_u = \partial T$.
\end{thmnonumber}
In Section \ref{existence}, we use this theorem to prove the existence of calibrations in a weak sense (see Definition \ref{calcur}) as a consequence of the Hahn-Banach theorem. The general idea of this proof follows closely Federer's approach to calibrations for minimal surfaces in \cite{RFC} and it suggests that, at least in dimension one, it would be possible to produce the analogue result and to extract an $L^\infty$ vector field playing the role of a calibration.\\
It is worth to notice that the convex decomposition formula presented in this paper relies on the one dimensional structure of the domain. In particular in Proposition \ref{coareasame} it is necessary that the singular points of an SBV function disconnect the domain; this is clearly peculiar of the dimension one, but it is likely that similar decomposition can be found in higher dimension and similar results could be obtained.\\
Moreover, even if all the proof of this paper are carried on for the functional \eqref{mum} the results can be extended with minor modifications to more general Mumford-Shah type functionals. We refer to Remark \ref{gener} for further details in this direction. 
\begin{ackn}
The author is warmly grateful to Professor Bernd Kirchheim for the useful discussions about this problem. The author would also like to thank Professor Giovanni Alberti for the valuable suggestions, Professor Domenico Mucci for the careful reading of the manuscript and the referee for the detailed review.
\end{ackn}

\section{Preliminaries}
Throughout the paper we consider $\Omega$ and $\Omega'$ to be open, bounded, regular sets of $\R^n$ such that $\Omega' \subset\subset \Omega$. Given $g\in L^\infty(\Omega)$ we define the Mumford-Shah functional as stated in the introduction 
\begin{equation}\label{mum2}
\mathfrak{F}(u) = \int_{\Omega} |\nabla u|^2 \, dx +  \mathcal{H}^{n-1}(S_u) + \int_{\Omega} |u-g|^2 \, dx
\end{equation}
and the homogeneous version
\begin{equation}\label{ho2}
F(u)= \int_{\Omega} |\nabla u|^2 \, dx +  \mathcal{H}^{n-1}(S_u),
\end{equation}   
where $u\in SBV(\Omega)$ and $S_u$ is the singular set of $u$. We refer to \cite{CCC} for the basic properties of BV and SBV functions and to \cite{SSO} for a comprehensive treatise on the Mumford-Shah functional. 
\vspace{2mm}\\
We deal with the following notions of minimizers:
\begin{defi}[Minimizer of $\mathfrak{F}$]
Given $g\in L^\infty(\Omega)$ we say that $u\in SBV(\Omega)$ is a minimizer of $\mathfrak{F}$ if $\mathfrak{F}(u) \leq \mathfrak{F}(v)$ for all $v\in SBV(\Omega)$. 
\end{defi}
\begin{defi}[Dirichlet minimizers]
We say that $u\in SBV(\Omega)$ is a Dirichlet minimizer of $F$ (resp. $\mathfrak{F}$) in $\Omega'$ if
\begin{displaymath}
F(u) \leq F(v) \quad \ \forall v\in SBV(\Omega)\quad \mbox{s.t.}\quad  u=v \quad \mbox{in} \quad \Omega\setminus \Omega'
\end{displaymath}
\begin{displaymath}
(resp. \,\ \mathfrak{F}(u) \leq \mathfrak{F}(v) \quad \ \forall v\in SBV(\Omega)\quad \mbox{s.t.} \quad  u=v \quad \mbox{in} \quad \Omega\setminus \Omega').
\end{displaymath}
\end{defi}
Proving that a function $u \in SBV(\Omega)$ is a Dirichlet minimizer in $\Omega'$ is not an easy question (in general); this is one of the main reasons why a notion of calibration resembling the one of minimal surfaces by Harvey and Lawson (\cite{CGE}) has turned out to be very useful. It was proposed by Alberti, Bouchittè and Dal Maso in \cite{TCM} and developed among the others in \cite{LCF} and \cite{LCFM}. In the next section we will give a brief introduction on this topic. 
\subsection{Calibration for the Mumford-Shah Functional}\label{sec}
Given $H:L^1(\Omega) \rightarrow \mathbb{R}$ let us define an abstract calibration in the following way:
\begin{defi}[Abstract calibration]\label{abstract}
Given $u\in L^1(\Omega)$, an abstract calibration for $u$ is a functional $G:L^1(\Omega) \rightarrow \mathbb{R}$ such that
\begin{equation}
(i) \quad H(u) = G(u), \quad \  (ii) \quad H(v) \geq G(v),\quad \  (iii) \quad G(u) = G(v) 
\end{equation} 
for all $v\in L^1(\Omega)$ such that $u=v$ in $\Omega\setminus \Omega'$.
\end{defi}
\begin{rmk}
If $G$ is a calibration for $u$, then $u$ is a Dirichlet minimizer in $\Omega'$ for $H$, indeed 
\begin{displaymath}
H(u) \overset{(i)}{=} G(u) \overset{(iii)}{=} G(v) \overset{(ii)}{\leq} H(v)
\end{displaymath}
for all $v \in L^1(\Omega)$  such that $u=v$ in $\Omega\setminus \Omega'$. 
\end{rmk}
In \cite{TCM} Alberti, Bouchitt\'e and Dal Maso introduced a stronger notion of calibration for the Mumford-Shah functional.
Given $v\in SBV(\Omega)$, we denote by $v^-(x)$ and $v^+(x)$ the lower and the upper traces of $v$. Moreover let $\Gamma_v$ be the extended graph of $v$ defined as
\begin{equation}\label{completegraph}
\Gamma_v = \{(x,t) \in \Omega \times \R: v^-(x) \leq t\leq v^+(x)\}.
\end{equation}
For standard theory on BV functions (\cite{CCC}) $\Gamma_v$ is rectifiable and then it admits a generalized normal that we are going to denote with $\nu_{\Gamma_v}$.\\
The calibration proposed in \cite{TCM} has the following form:
\begin{displaymath}
G(v) = \int_{\Gamma_v} \langle \phi , \nu_{\Gamma_v}\rangle \, d\mathcal{H}^n,
\end{displaymath}
where $\phi : \Omega \times \mathbb{R} \rightarrow \mathbb{R}^{n+1}$ is a vector field to be determined. The regularity asked on $\phi$ is the least that guarantees the validity of a divergence theorem on $\Omega \times \R$. To be more precise we refer to \cite{TCM} and for reader convenience we propose the definition of $\emph{approximately regular}$ vector field:
\begin{defi}[Approximately regular vector field]
Given $A\subset \R^{n+1}$, a vectorfield $\phi: A \rightarrow \R^{n+1}$ is approximately regular if it is bounded and for every Lipschitz hypersurface $M$ in $\R^{n+1}$ there holds
\begin{equation}
\lim_{r\rightarrow 0}\mint_{B_r(x_0)\cap A} |(\phi(x)-\phi(x_0))\cdot \nu_M(x_0)|\, dx = 0 
\end{equation}
for $\mathcal{H}^{n}$-a.e. $x_0\in M\cap A$.
\end{defi}  
Comparing the functional $G$ with $F$, it is possible to find sufficient conditions on $\phi$ such that $G$ satisfies properties (i), (ii) and (iii) with respect to $F$ for a given $u\in SBV(\Omega)$. Then the vector field satisfying these properties is called \emph{calibration} for $u$.
\begin{defi}[Calibration for the Mumford-Shah Functional, \cite{TCM}]
Let $\Omega \subset \mathbb{R}^n$ be open and bounded and $u\in SBV(\Omega)$. Given $\phi = (\phi^x,\phi^t) : \Omega \times \mathbb{R} \rightarrow \mathbb{R}^{n+1}$ an approximately regular vector field, we say that it is a calibration for $u$ if it is divergence free and 
\begin{itemize}
\item[a)] $\phi^t(x,t) \geq \dfrac{|\phi^x(x,t)|^2}{4} \quad \mbox{for } \mathscr{L}^n$-a.e. $x\in \Omega$ and for all $t\in \mathbb{R}$,
\item[b)] $\left|\displaystyle\int_{t_1}^{t_2} \phi^x(x,t) \, dt\right| \leq 1 \quad \mbox{for } \mathcal{H}^{n-1}$-a.e. $x\in \Omega$ and for all $t_1,t_2\in \mathbb{R}$, \\
\item[c)]$\phi^x(x,u(x) ) = 2\nabla u(x), \qquad \phi^t(x,u(x)) = |\nabla u(x)|^2 \quad \mbox{for } \mathscr{L}^n$-a.e. $x\in \Omega$, \\
\item[d)] $\displaystyle\int_{u^-(x)}^{u^+(x)} \phi^x(x,t) \, dt = \nu_u(x) \quad \mbox{for } \mathcal{H}^{n-1}$-a.e. $x\in S_u$,
\end{itemize}
where $\nu_u$ is the approximate normal of $S_u$.
\end{defi} 
As properties $(a),(b),(c),(d)$ imply $(i),(ii)$ and $(iii)$ for $G$ we have the following theorem:
\begin{thm}[\cite{TCM}]
Given $u \in SBV (\Omega)$, suppose that there exists $\phi: \Omega \times \mathbb{R} \rightarrow \mathbb{R}^{n+1}$ a calibration for $u$. Then $u$ is a Dirichlet minimizer in $\Omega'$ of the homogeneous Mumford-Shah functional \eqref{ho2}.
\end{thm}
In an analogous way a similar notion can be introduced in order to study minimizers of $\mathfrak{F}$. It is enough to replace conditions $(a)$ and $(c)$ with
\begin{itemize}
\item[a')]
$\phi^t(x,t) \geq \dfrac{|\phi^x(x,t)|^2}{4} - (t-g)^2 \quad \mbox{for } \mathscr{L}^n$-a.e. $x\in \Omega$ and for all $t\in \mathbb{R}$,\\
\item[c')] $\phi^x(x,u(x) ) = 2\nabla u(x) , \quad \phi^t(x,u(x)) = |\nabla u(x)|^2 - (u-g)^2\quad \mbox{for } \mathscr{L}^n$-a.e. $x\in \Omega$.
\end{itemize}
\begin{thm}[\cite{TCM}]
Given $u \in SBV (\Omega)$, suppose that there exists $\phi: \Omega \times \mathbb{R} \rightarrow \mathbb{R}^{n+1}$ a calibration for $u$ with $(a)$ and $(c)$ replaced with $(a')$ and $(c')$. Then $u$ is a Dirichlet minimizer in $\Omega'$ of the Mumford-Shah functional \eqref{mum2}.
\end{thm}
As a consequence, in \cite{TCM}, the authors proposed the following alternative formulation of the Mumford-Shah functional 
\begin{equation}\label{formula1}
F(u)= \max_{\phi\in K} \int_{\Gamma_u} \langle \phi, \nu_{\Gamma_u}\rangle \, d \mathcal{H}^{n}= \max_{\phi\in K} \int_{\Omega\times \mathbb{R}} \langle \phi, D\textbf{1}_{\{u>t\}}\rangle , 
\end{equation}
\begin{equation}\label{formula2}
\mathfrak{F}(u)= \max_{\phi\in K'} \int_{\Gamma_u} \langle \phi , \nu_{\Gamma_u} \rangle \, d \mathcal{H}^{n}= \max_{\phi\in K'} \int_{\Omega\times \mathbb{R}} \langle \phi ,D\textbf{1}_{\{u>t\}}\rangle, 
\end{equation}
where 
\begin{equation}\label{kappa1}
K = \{\phi:\Omega\times\mathbb{R}\rightarrow \mathbb{R}^{n+1},\ Borel: (a) \mbox{ and }(b)\mbox{ hold pointwise} \}
\end{equation}
and
\begin{equation}\label{kappa2}
K' = \{\phi: \Omega\times\mathbb{R}\rightarrow \mathbb{R}^{n+1},\ Borel: (a') \mbox{ and }(b)\mbox{ hold pointwise}\}.
\end{equation}
\begin{rmk}
The previous representation formula is the starting point for the proof of existence of calibration in dimension one, due to Chambolle \cite{CRF}. In particular one can introduce the following convex functional also called $\emph{lift}$ of $F$
\begin{displaymath}
\mathcal{F}_K(w)= \sup_{\phi\in K \cap C_0(\Omega \times \R, \R^{n+1})} \int_{\Omega\times\mathbb{R}} \phi \cdot Dw, 
\end{displaymath} 
with $w:I\times \R \rightarrow [0,1]$ decreasing in the second variable and of bounded variation. In \cite{CRF} Chambolle proves that if $u\in SBV(I)$ is a minimizer of the Mumford-Shah functional then $\textbf{1}_{\{u(x)>t\}}$ is a minimizer of $\mathcal{F}_K$. Then by Hahn-Banach theorem it is possible to prove the existence of calibrations in a weak asymptotic sense. 
\end{rmk}
\begin{rmk}\label{coareafalse}
It is interesting to notice that one can prove the same result in higher dimension if $\mathcal{F}_K$ satisfies a generalized coarea formula of the form
\begin{equation}\label{chacoa}
\mathcal{F}_K(w) = \int_0^1 \mathcal{F}_K(\textbf{1}_{\{w(x,t)>s\}})\, ds.
\end{equation}
Unfortunately this is false even in dimension one. Indeed it is enough to consider
\begin{displaymath}
u_1(x)=\left\{
\begin{array}{ll}
0 & \mbox{ if } x\leq 1/2\\
x & \mbox{ if } x> 1/2,
\end{array}
\right.
\qquad
u_2(x)=\left\{
\begin{array}{ll}
x & \mbox{ if } x\leq 1/2\\
1 & \mbox{ if } x> 1/2
\end{array}
\right.
\end{displaymath}
and 
$w(x,t) = (1/2)\textbf{1}_{\{u_1(x)> t\}} + (1/2)\textbf{1}_{\{u_2(x)> t\}}$ to see that formula \eqref{chacoa} does not hold.
\end{rmk}

\subsection{A lifting of the Mumford-Shah functional in the space of rectifiable currents}
In this section we introduce a lifted functional that takes values in $\mathcal{R}_{n}(\Omega \times \R)$ the $n$-dimensional rectifiable currents with real multiplicity. We briefly recall the basic theory of currents and we refer the reader to \cite{CCC} for a more detailed overview.\\
Let $U$ be an open subset of $\R^N$. A $k$-dimensional current on $U$ is a linear continuous (see \cite{CCC}) functional on the space of $k$-forms $\Lambda^k(U)$ with coefficients in $C_c^\infty(U)$. \\
In particular we define the space $\mathcal{R}_{k}(U)$ of $k$-dimensional rectifiable currents with real multiplicity as the triple $(\mathcal{M},\theta,\xi)$ where $\mathcal{M} \subset U$ is a $k$-rectifiable set, $\theta : \mathcal{M} \rightarrow \R_+$ is a function called multiplicity and $\xi$ is a map that associates to $\mathcal{H}^n$-a.e. $x$ in $\mathcal{M}$ a unit, simple $k$-vector orienting $\mathcal{M}$.
We define the current $(\mathcal{M},\theta,\xi)$ by its action on a $k$-diffential form $\omega \in \Lambda^k(U)$ in the following way:
\begin{displaymath}
(\mathcal{M},\theta,\xi)(\omega) = \int_{\mathcal{M}} \langle \omega, \xi\rangle \theta \, d\mathcal{H}^k,
\end{displaymath}
where $\langle\cdot , \cdot \rangle$ denote the duality product between vectors and covectors.
Moreover given $T=(\mathcal{M},\theta,\xi)$ we define the total variation measure associated to $T$ as
\begin{displaymath}
\|T\| (A) = \int_{\mathcal{M}\cap A} \theta \, d\mathcal{H}^k
\end{displaymath} 
for every $A \subset U$ measurable. We call $\|T\|(U) = M(T)$ the mass of $T$.\\
We define the restriction of a rectifiable current $T=(\mathcal{M},\theta,\xi)$ on a measurable set as 
\begin{displaymath}
T\res A(\omega) = \int_{\mathcal{M}\cap A} \langle \omega, \xi\rangle \theta \, d\mathcal{H}^k
\end{displaymath}
for every $A \subset U$ measurable. In addition given $\alpha\in \Lambda^h(U)$ with $h\leq k$, we define the restriction of $T \in \mathcal{R}_k(U)$ to $\alpha$ as the $(k-h)$-dimensional current $T\res \alpha$ defined as
\begin{displaymath}
T\res \alpha(\omega) = T(\alpha \wedge \omega) 
\end{displaymath} 
for every $\omega\in \Lambda^{k-h}(U)$.\\
Moreover, given $E$ a $k$-rectifiable set in $U$ we will denote by $[\![E]\!]$ the $k$-dimensional rectifiable current induced by $E$, that is defined as
\begin{displaymath}
[\![E]\!](\omega) := \int_{E} \langle\omega, \xi_{E}\rangle \, d\mathcal{H}^k,
\end{displaymath}
for $\omega \in \Lambda^k(U)$, where $\xi_{E}$ is the unit simple $k$-vector orienting $E$.
Therefore, as a consequence of \eqref{completegraph}, we can define $[\![\Gamma_u]\!]$, the $n$-rectifiable current associated to the complete graph of $u\in SBV(\R^n)$. From now on, with a little abuse of notation, we will denote it by $\Gamma_u$, instead of $[\![\Gamma_u]\!]$ (it will be clear by the context if we are dealing with the rectifiable set or with the current associated to it).
\vspace*{2mm}\\
We introduce the lifting of the Mumford-Shah functional on the space of rectifiable currents for the functionals $\mathfrak{F}$ and $F$.
\begin{defi}[Lifting to the space of rectifiable current]
Given $T=(\mathcal{M},\theta,\xi) \in \mathcal{R}_{n}(\Omega \times \R)$ we define
\begin{equation}\label{func}
G_K(T) := \sup_{\phi\in K} \int_{\mathcal{M}} \langle \phi , \star (-\xi) \rangle d\|T\| =\sup_{\phi\in K} \int_{\mathcal{M}} \theta\langle \phi , \nu_T \rangle d\mathcal{H}^{n}  
\end{equation}
and
\begin{equation}
G_{K'}(T) := \sup_{\phi\in K'} \int_{\mathcal{M}} \langle \phi , \star (-\xi) \rangle d\|T\| = \sup_{\phi\in K} \int_{\mathcal{M}} \theta\langle \phi , \nu_T \rangle d\mathcal{H}^{n} 
\end{equation}
where $\nu_T := - (\star \xi)$, $\star$ is the Hodge star and $K$ and $K'$ are defined as in \eqref{kappa1} and in \eqref{kappa2}.
\end{defi}
\begin{prop}\label{convexlsc}
The functionals $G_K$ and $G_{K'}$ satisfy the following properties:
\begin{itemize}
\item [$(i)$] They are convex on $\mathcal{R}_{n}(\Omega \times \R).$
\item [$(ii)$] They are lower semicontinous with respect to the mass bounded convergence.
\item [$(iii)$] Given $v\in SBV(\Omega)$, $G_K(\Gamma_v) = F(v)$ and $G_{K'}(\Gamma_v) = \mathfrak{F}(v)$.
\end{itemize}
\end{prop}
\begin{dimo}
Statement $(i)$ follows from the definition and $(iii)$ is a consequence of the representation formulas \eqref{formula1} and \eqref{formula2}.
Moreover $(ii)$ can be proved with an easy modification of the argument in \cite{CCC} sec. $3.3.1$.
\end{dimo}

\section{A convex decomposition formula for the Mumford-Shah functional in dimension one}\label{mais}
We restrict our analysis to the case $n=1$. We also assume $\Omega = I$ and $\Omega' = I'$ to be open and bounded intervals such that $I' \subset \subset I$ and we consider the Mumford-Shah functional in its general form
\begin{equation}\label{generic}
F(u) := \int_I |u'(x)|^2\, dx + \beta\int_I|u-g|^2\, dx + \alpha\mathcal{H}^0(S_u), 
\end{equation} 
where $\alpha > 0$, $\beta \geq 0$, $g\in L^\infty(I)$ and $u\in SBV(I)$. Notice that when $\beta = 0$ and $\alpha = 1$, $F$ is the homogeneous version of the Mumford-Shah functional as defined in \eqref{ho2}. From now on we will denote by $u^l(x)$ (resp. $u^r(x)$) the left (resp. right) trace of $u$ in a point $x$.

\begin{rmk}\label{gener}
Even if we restrict our attention to \eqref{generic} it is important to remark that the results of this section and of the following one hold for a more general class of functionals with minor modification of the proofs. Functionals of the form
\begin{displaymath}
W(u) = \int_I f(u'(x),u(x),x)\, dx + \sum_{x\in S_u} \psi(x,u^l(x),u^r(x))  
\end{displaymath} 
with suitable hypothesis on $f$ and $\psi$ necessary to ensure the lower semicontinuity of $W$ and the existence of minimizers can be treated by this theory. We refer to \cite{ETF} for the precise assumptions and we stress the fact that in our setting $f$ need not to be assumed more regular as in \cite{CRF}. For example in the case of the Mumford-Shah functional $g$ can be taken in $L^\infty$ without affecting the proof, while in \cite{CRF} the function $g$ needs to have a l.s.c. and a u.s.c. representatives in $L^\infty$.
\end{rmk}
If we consider the functional $F$ as defined in \eqref{generic}, its convex lift defined in \eqref{func} on $\mathcal{R}_1(I\times\R)$ reads
\begin{equation}\label{operator}
G(T) =  \sup_{\phi\in K} \int_{\mathcal{M}} \theta\langle \phi , \nu_T \rangle d\mathcal{H}^1 
\end{equation}
for every $T=(\mathcal{M},\theta,\xi)$. \\
In particular $K$ is the set of $\phi :I\times \R\rightarrow \R^2$, Borel, such that
\begin{itemize}
\item[I)] $\phi^t(x,t) \geq \dfrac{|\phi^x(x,t)|^2}{4} - \beta(t-g)^2 \quad \mbox{for all }x\in I$ and for all $t\in \mathbb{R}$,\\
\item[II)] $\left|\displaystyle\int_{t_1}^{t_2} \phi^x(x,t) \, dt\right| \leq \alpha \quad \mbox{for all } x\in I$ and for all $t_1,t_2\in \mathbb{R}$.
\end{itemize}
We are going to consider as the domain of $G$ the cone $C \subset \mathcal{R}_1(I \times \R)$ made by finite linear combinations of SBV graphs:
\begin{equation}\label{cone}
C:= \left\{T= \sum_{i=1}^{k} \lambda_i \Gamma_{u_i} : k\in \mathbb{N}, \lambda_i \in \mathbb{R}_+, u_i \in SBV(I)\right\}.
\end{equation}
For every $T\in C$ we will assume implicitly that, being a rectifiable current, it is defined by a triple $T=(\mathcal{M},\theta,\xi)$.

\subsection{Simplifying the cone $C$}
From the definition of the cone $C$ in \eqref{cone} one easily notices that for every current $T\in C$ there exists different combinations of SBV graphs $\{u_i\}$ that represent it. In particular there are some configurations we would like to avoid and this subsection is devoted to make this simplifications for $C$.  
\begin{defi}
Given $\{u_i\}_{i=1\ldots k} \subset SBV(I)$. We say that the family $\{u_i\}_{i=1\ldots k}$ has \emph{cancellation on the jumps} if there exists $l_1,l_2$ and $x_0 \in S_{u_{l_1}} \cap S_{u_{l_2}}$ such that
\begin{displaymath}
u_{l_1}^l(x_0) \leq  u_{l_2}^r(x_0) < u_{l_1}^r(x_0) \leq u_{l_2}^l (x_0) \quad \mbox{or}
\end{displaymath} 
\begin{displaymath}
u_{l_1}^r(x_0) \leq  u_{l_2}^l(x_0) < u_{l_1}^l(x_0) \leq u_{l_2}^r (x_0).
\end{displaymath} 
\end{defi} 
We need a lemma that ensures that we can rearrange the graphs in order not to have this cancellation.
\begin{lemma}\label{nocanc}
Given $T =\sum_{i=1}^k \lambda_i \Gamma_{u_i} \in C$ there exists $l \in \mathbb{N}$, $w_i \in SBV(I)$ and $\mu_i \in \R^+$ for $i=1\ldots l$ such that $T = \sum_{i=1}^l \mu_i \Gamma_{w_i}$ and
there is no cancellation on the jumps. 
\end{lemma}
\begin{dimo}
Given $T= \sum_{i=1}^{k} \lambda_i \Gamma_{u_i}$ let us suppose that we have cancellation between $\Gamma_{u_1}$ and $\Gamma_{u_2}$ in $A\subset S_{u_1} \cap S_{u_2}$ and $\lambda_1 \geq \lambda_2$ (without loss of generality).
As $A$ is countable we will denote it by the sequence $\{x_1, x_2, \ldots \}$ possibly infinite.
Given $I=(a,b)$ consider the new sequence $\{a=x_0, x_1, x_2, \ldots \}$ and define two SBV functions in the following way: 
\begin{displaymath}
w_1(x) = \left\{
\begin{array}{ll}
u_1(x) & \mbox{ for } x_{i-1}<x\leq x_{i},\ i\geq 1 \mbox{ and odd}\\
u_2(x) & \mbox{ for } x_{i-1}<x\leq x_{i},\ i\geq 1 \mbox{ and even}
\end{array}
\right.
\end{displaymath}
and 
\begin{displaymath}
w_2(x) = \left\{
\begin{array}{ll}
u_2(x) & \mbox{ for } x_{i-1}<x\leq x_{i},\ i\geq 1 \mbox{ and odd }\\
u_1(x) & \mbox{ for } x_{i-1}<x\leq x_{i},\ i\geq 1 \mbox{ and even}.
\end{array}
\right.
\end{displaymath}
Then we have that $\lambda_2 \Gamma_{w_1} + \lambda_2\Gamma_{w_2} + (\lambda_1 - \lambda_2) \Gamma_{u_1} =  \lambda_1\Gamma_{u_1}+ \lambda_2\Gamma_{u_2}$. Hence we produce a decomposition of $\lambda_1\Gamma_{u_1}+ \lambda_2\Gamma_{u_2}$ that has no cancellation on the jumps. It is easy to check that one can repeat this operation for any pair of graphs that has cancellation on jumps and that this procedure ends in a finite number of steps.
\end{dimo}

\subsection{Properties of the regular part of $G(T)$}
\begin{defi}[Regular part and singular part of $T$]
We define the singular part of $T = \sum_{i=1}^k\lambda_i \Gamma_{u_i} \in C$ as 
\begin{equation}
S_T := \bigcup_{i=1}^k S_{u_i}
\end{equation}
and the regular part as $R_T := I \setminus S_T$.
\end{defi}
\begin{rmk}
One can easily notice that if we assume that the graphs do not have cancellation according to Lemma \ref{nocanc}, $S_T$ is well defined, so it does not depend on the representation of $T$.
\end{rmk}
Given a measurable set $A\subset I$ we define the localized version of $G$ as
\begin{displaymath}
G(T,A) := \sup_{\phi\in K} \int_{\mathcal{M}\cap (A\times \R)} \langle \phi , \nu_T \rangle d\|T\|. 
\end{displaymath}
\begin{rmk}
It is clear that given $A_1$, $A_2$ disjoint measurable sets we have
\begin{displaymath}
G(T,A_1 \cup A_2) = G(T,A_1) + G(T,A_2)
\end{displaymath}
so in particular
\begin{equation}
G(T) = G(T,S_T) + G(T,R_T).
\end{equation}
Moreover when one computes the localized functional, it is possible to restrict the set $K$ accordingly:
\begin{displaymath}
G(T,A) = \sup_{\phi\in K_A} \int_{\mathcal{M}\cap (A\times \R)} \langle \phi , \nu_T \rangle d\|T\|,
\end{displaymath}
where $K_A$ is the set of $\phi :I\times \R\rightarrow \R$, Borel, such that
\begin{itemize}
\item
$\phi^t(x,t) \geq \dfrac{|\phi^x(x,t)|^2}{4} - \beta(t-g)^2 \quad \forall x\in A$ and $\forall t\in \mathbb{R}$,\\
\item $\left|\displaystyle\int_{t_1}^{t_2} \phi^x(x,t) \, dt\right| \leq \alpha \quad \mbox{for every } x\in A$ and for all $t_1,t_2\in \mathbb{R}$. 
\end{itemize}
\end{rmk}
We are presenting a proposition that allows us to split $G(T,R_T)$ as the sum of $\lambda_i G(\Gamma_{u_i},R_T)$.
\begin{prop}\label{decomposition}
Given $T = \sum_{i=1}^k \lambda_i \Gamma_{u_i} \in C$, then
\begin{equation}\label{regular}
G(T, R_T) = \sum_{i=1}^k \lambda_i G(\Gamma_{u_i}, R_T)=\sum_{i=1}^k \lambda_i\left(\alpha\int_I (u_i')^2 \, dx + \beta\int_I |u_i - g|^2 \,dx\right).  
\end{equation}
\end{prop}
In order to give a proof of this fact we need some preliminary lemmas.
\begin{lemma}\label{first}
Given $T = \sum_{i=1}^k \lambda_i \Gamma_{u_i} \in C$ let $A \subset I$ be a measurable set such that $A \cap S_T = \emptyset$ and $\mathcal{H}^{1}(\Gamma_{u_i} \cap \Gamma_{u_j} \cap (A\times \R)) = 0$ for every $i\neq j$. Then 
\begin{displaymath}
G(T,A) = \sum_i \lambda_i G(\Gamma_{u_i},A).
\end{displaymath}
\end{lemma}
\begin{dimo}
By induction it is enough to show that given, $T_1 = \sum_{i=1}^{k-1}\lambda_i \Gamma_{u_i}$ and $T_2= \lambda_k \Gamma_{u_k}$ one has
\begin{displaymath}
G(T_1 + T_2 ,A) = G(T_1,A) + G(T_2,A). 
\end{displaymath} 
Fix $\varepsilon>0$. For $i=1,2$ there exist $\phi_i \in K_A$ such that
\begin{displaymath}
\int_{\mathcal{M}_i \cap (A\times \R)} \langle \phi_i , \nu_{T_i} \rangle \, d\|T_i\| \geq  G(T_i, A) - \varepsilon,
\end{displaymath}
where  $T_i = (\mathcal{M}_i, \theta_{T_i} ,\xi_{T_i})$ and $\nu_{T_i} = -(\star \xi_{T_i})$.
Define then the following vector field
\begin{displaymath}
\tilde \phi = \left\{\begin{array}{ll}
\phi_1 & (x,t) \in \mathcal{M}_1 \setminus \mathcal{M}_2\\
\phi_2 & (x,t) \in \mathcal{M}_2\setminus \mathcal{M}_1 \\
 0 & \mbox{ otherwise.}
\end{array}
\right.
\end{displaymath}
Let us prove that $\tilde \phi \in K_A$. \\
For every $x \in A$ we have that $x\notin S_{T}$ by hypothesis, so that (II) is satisfied and (I) is trivial by definition.
Moreover, as $\mathcal{H}^1(\mathcal{M}_1\cap \mathcal{M}_2 \cap (A\times \R)) = 0$, one has 
\begin{eqnarray*}
\int_{(\M_1 \cup \M_2)\cap (A\times \R)}  \langle \tilde \phi ,\nu_T \rangle \, d\mathcal{H}^1 = \int_{\M_1 \cap (A\times \R)} \langle \phi_1 , \nu_{T_1} \rangle  \, d\mathcal{H}^1 + \int_{\M_2\cap (A\times \R)} \langle \phi_2 , \nu_{T_2} \rangle  \, d\mathcal{H}^1.
\end{eqnarray*}
So
\begin{eqnarray*}
G(T_1,A) + G(T_2,A) &\leq &\int_{\M_1 \cap (A\times \R)} \langle \phi_1 , \nu_{T_1} \rangle \theta_{T_1}  \, d\mathcal{H}^1 + \int_{\M_2\cap (A\times \R)} \langle \phi_2 , \nu_{T_2} \rangle \theta_{T_2} \, d\mathcal{H}^1 + 2\varepsilon \\
&\leq & G(T_1 + T_2 ,A) + 2\varepsilon. 
\end{eqnarray*}
Sending $\varepsilon$ to zero we obtain the first inequality.
The opposite one comes directly from the convexity of $G$. 
\end{dimo}

\begin{figure}
\centering
\subfigure{
\begin{tikzpicture}[>=stealth, scale=0.9]
\draw [->](-1,0) -- (6,0) node [scale=0.7][below]{x};
\draw [->](0,-1) -- (0,6) node [scale=0.7][left]{t};
\draw[densely dashed] (5,0) -- (5,6);
\draw  (0,0) [thick][red]
to[out=10,in=180, looseness=1] (2,2) 
to[out=90,in=-90, looseness=1] (2,2.5) 
to[out=0,in=160, looseness=0.8] (5,5) node[black, scale=0.7] [below=18pt, left=4pt] {$\Gamma_{u_1}$} ;

\draw  (0,0)[thick][blue]
to[out=40,in=180, looseness=1] (1.5,3.5)
to[out=90,in=-90, looseness=1] (1.5,4)
to[out=0,in=160, looseness=0.8] (5,5)node[black, scale=0.7] [above=3pt,left=35pt]{$\Gamma_{u_2}$};

\draw  (2,0) node[black, scale=0.7] [below=0.5pt]{$S_{u_1}$} [dashed] -- (2,2) ;
\draw  (1.5,0) node[black, scale=0.7] [below=0.5pt]{$S_{u_2}$} [dashed] -- (1.5,3.5) ;

\draw  (0,2) node[black, scale=0.7] [left=0.5pt]{$u_1^l$} [dashed] -- (2,2) ;
\draw  (0,2.5) node[black, scale=0.7] [left=0.5pt]{$u_1^r$} [dashed] -- (2,2.5) ;
\draw  (0,3.5) node[black, scale=0.7] [left=0.5pt]{$u_2^l$} [dashed] -- (1.5,3.5) ;

\draw  (0,4) node[black, scale=0.7] [left=0.5pt]{$u_2^r$} [dashed] -- (1.5,4) ;

\end{tikzpicture}
}\qquad\subfigure{\begin{tikzpicture}[>=stealth, scale=0.9]
\draw [->](-1,0) -- (6,0) node [scale=0.7][below]{x};
\draw [->](0,-1) -- (0,6) node [scale=0.7][left]{t};
\draw[densely dashed] (5,0) -- (5,6);
\draw  (1,1) [thick][red]
to[out=10,in=180, looseness=1] (2,2) 
to[out=90,in=-90, looseness=1] (2,2.5) 
to[out=0,in=160, looseness=0.8] (5,5) node[black, scale=0.7] [below=14pt, left=4pt] {$\Gamma_{u_1}$} ;

\draw  (1,1)[thick][blue]
to[out=65,in=250, looseness=1] (1.5,3.5)
to[out=90,in=-90, looseness=1] (1.5,4)
to[out=0,in=160, looseness=0.8] (5,5)node[black, scale=0.7] [above=6pt,left=35pt]{$\Gamma_{u_2}$};

\draw  (0,0)[thick][green]
to[out=40,in=180, looseness=1] (1,1)node[black, scale=0.7] [left=7pt,below=5pt]{$\Gamma_{u_1} = \Gamma_{u_2}$};

\draw  (2,0) node[black, scale=0.7] [below=0.5pt]{$S_{u_1}$} [dashed] -- (2,2) ;
\draw  (1.5,0) node[black, scale=0.7] [below=0.5pt]{$S_{u_2}$} [dashed] -- (1.5,3.5) ;

\draw  (0,2) node[black, scale=0.7] [left=0.5pt]{$u_1^l$} [dashed] -- (2,2) ;
\draw  (0,2.5) node[black, scale=0.7] [left=0.5pt]{$u_1^r$} [dashed] -- (2,2.5) ;
\draw  (0,3.5) node[black, scale=0.7] [left=0.5pt]{$u_2^l$} [dashed] -- (1.5,3.5) ;

\draw  (0,4) node[black, scale=0.7] [left=0.5pt]{$u_2^r$} [dashed] -- (1.5,4) ;

\end{tikzpicture}
}
\caption{Configuration in Lemma \ref{first} and \ref{overlapping}}
\label{firstover}
\end{figure}
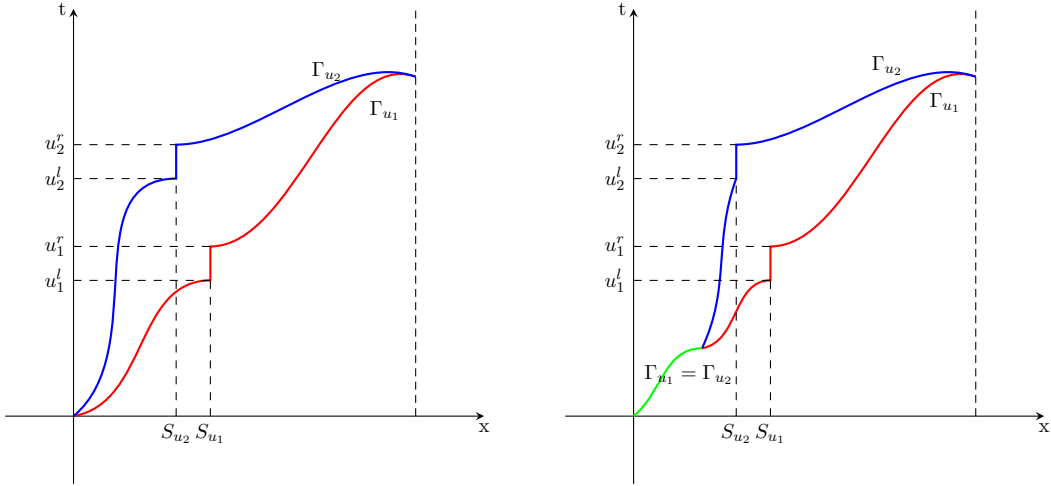

\begin{lemma}\label{overlapping}
Given $T = \sum_{i=1}^k \lambda_i \Gamma_{u_i} \in C$ let $A\subset I$ be a measurable set such that $A \cap S_T = \emptyset$. Then 
\begin{displaymath}
G(T,A) = \sum_i \lambda_i G(\Gamma_{u_i},A).
\end{displaymath}
\end{lemma}
\begin{dimo}
Given $T\in C$, let $J$ be a set of indexes. Denote by $\Gamma = \bigcap_{i\in J} \Gamma_{u_i}$ an intersection of graphs and let $\theta = \sum_{i\in J} \lambda_i$ be the multiplicity on $\Gamma$. So
\begin{eqnarray*}
 \sup_{\phi\in K} \int_{\Gamma \cap (A\times \R)} \langle \phi,\nu_T\rangle \, d\|T\|&=& \sup_{\phi\in K} \int_{\Gamma\cap (A\times \R)}  \theta\langle \phi,\nu_T\rangle \, d\mathcal{H}^1 = \sup_{\phi\in K} \int_{\Gamma \cap (A\times \R)} \sum_{i\in J} \lambda_i\langle \phi,\nu_T\rangle \, d\mathcal{H}^1 \\
&=& \sum_{i\in J} \lambda_i\sup_{\phi\in K} \int_{\Gamma \cap (A\times \R)}  \langle \phi,\nu_T\rangle \, d\mathcal{H}^1.
\end{eqnarray*}
Clearly this can be repeated for every intersection of an arbitrary number of graphs. Combining this result with Lemma \ref{first} we have the thesis.
\end{dimo}
\emph{Proof of Proposition \ref{decomposition}}
\vspace{2mm}\\
Proposition \ref{decomposition} is a direct consequence of Lemma \ref{overlapping} choosing $A = R_T$ and the second equality in \eqref{regular} follows from Proposition \ref{convexlsc}.
\begin{flushright}
$\Box$
\end{flushright} 
\subsection{Properties of the singular part of $G(T)$}
In this section we are going to study the properties of $\mathcal{G}(T):= G(T,S_T)$. \\
Given $T= \sum_{i=1}^k \lambda_i \Gamma_{u_i}\in C$ and calling $\nu_T= ((\nu_T)^x, (\nu_T)^t)$, by \eqref{operator} we have 
\begin{displaymath}
\mathcal{G}(T) = \sup_{\phi \in K}\int_{\mathcal{M} \cap (S_T \times \R)} \theta \phi^x (\nu_{T})^x \, d\mathcal{H}^1 
\end{displaymath}
and it is easy to see that
\begin{displaymath}
(\nu_{T})^x(x,t) = \left\{\begin{array}{ll}
+1 & (x,t) \in S_{u_i} \times (u_i^l, u_i^r)\\ 
-1 & (x,t) \in S_{u_i} \times (u_i^r, u_i^l).
\end{array}
\right.
\end{displaymath}
Hence
\begin{displaymath}
\mathcal{G}(T) = \sup_{\phi \in K} \sum_{i=1}^k\int_{S_{u_i} \times (u_i^l, u_i^r)} \theta\phi^x\, d\mathcal{H}^1. 
\end{displaymath}
From now on we will work with linear combinations of graphs with the same multiplicity. We will see later the reason why we can reduce to this situation. We want to prove that, given $T= \sum_i \Gamma_{u_i}$, $\mathcal{G}(T)$ can be written as the sum of $\mathcal{G}(\Gamma_{u_i})$ in all the configurations in which there is non-adjacency of the jumps of the graphs. 
\begin{thm}\label{split}
Consider $T \in C$ such that $T = \sum_{i=1}^k \Gamma_{u_i}$. Suppose that for every $i,j = 1\ldots k$
\begin{displaymath}
\{x\in S_{u_i} \cap S_{u_{j}} : u^r_{i} (x) = u^l_{j}(x)\} = \emptyset. 
\end{displaymath}
Then
\begin{displaymath}
\mathcal{G}\left(\sum_{i=1}^k \Gamma_{u_i}\right) = \sum_{i=1}^k \mathcal{G}\left(\Gamma_{u_i}\right).
\end{displaymath}
\end{thm}
\begin{rmk}
Notice that without loss of generality we can prove the previous statement restricting the functional $\mathcal{G}$ to every $x\in S_T$. So the lemmas needed to prove Theorem \ref{split} will be stated for a fixed point $x\in S_T$.
\end{rmk}
For sake of clarity we propose two lemmas (Lemma \ref{fillingone} and \ref{fillingtwo}) that deals with a simple situation that is enough to explain the general strategy (See Figure \ref{splitfig}). Then, in Proposition \ref{summ1} and \ref{summ2}, we generalize this procedure and finally we prove the theorem.
\begin{lemma}\label{fillingone}
Consider $T= \sum_{i=1}^k \Gamma_{u_i} \in C$ such that $u_i$ are ordered in an increasing way. Fix $x\in S_T$ and suppose that we have $u_i^l(x) \leq u_i^r(x)$ for every $i=1 \ldots k$. Suppose in addition that
\begin{displaymath}
u_i^r(x) < u_j^l(x)\quad \mbox{ for every }  \quad i<j.
\end{displaymath}
Then
\begin{displaymath}
\mathcal{G}(T,\{x\}) = \sum_{i=1}^k\mathcal{G}(\Gamma_{u_i},\{x\})= \alpha |\{i: x\in S_{u_i}\}|. 
\end{displaymath}
In addition the maximum is achieved and letting $\phi_{T}$ be the vector field realizing the maximum for $T$
\begin{displaymath}
\phi^x_T(x,t) = \alpha/(u_i^r - u_i^l) \quad \mbox{ for every } \quad t\in (u_i^l,u_i^r)
\end{displaymath}
for every $i=1\ldots k$ such that $x\in S_{u_i}$.
\end{lemma}
\begin{dimo}
First of all notice that it is not restrictive to assume that $x	\in S_{u_i}$ for every $i=1,\ldots,k$. By induction it is enough to prove that for $T=T_1 + T_2$ where $T_1 = \sum_{i=1}^{k-1} \Gamma_{u_i}$ and $T_2 = \Gamma_{u_k}$ one has 
\begin{displaymath}
\mathcal{G}(T_1 + T_2,\{x\}) = \mathcal{G}(T_1,\{x\}) + \mathcal{G}(T_2,\{x\}) 
\end{displaymath}
and
\begin{displaymath}
\phi^x_{T}(x,t) = \alpha/(u_k^r - u_k^l)  \quad \mbox{ for every }  \quad t\in (u_k^l,u_k^r).
\end{displaymath}
(We suppose $x\in S_{u_k}$ because if not, there is nothing to prove).\\
For the inductive hypothesis we have that for all $i=1\ldots k-1$
\begin{displaymath}
\phi^x_{T_1}(x,t) = \alpha/(u_i^r - u_i^l) \quad \mbox{for every }  \quad t\in (u_i^l,u_i^r).
\end{displaymath}
For the general theory of calibration we have that, calling $\phi_{T_2}$ the vector field realizing the maximum in $\mathcal{G}(T_2,\{x\})$,
\begin{displaymath}
\phi^x_{T_2}(x,t) = \alpha/(u_k^r - u_k^l)  \quad \mbox{for every }  \quad t\in (u_k^l,u_k^r),
\end{displaymath}
because
\begin{displaymath}
\int_{u_k^-}^{u_k^+} \phi^x_{T_2}(x) = \alpha \quad \mbox{for every } x \in S_{u_k}.
\end{displaymath} 
Define the following vector field on $\{x\} \times \R$:
\begin{displaymath}
\tilde \phi = \left\{\begin{array}{ll}
\phi_{T_1} & (x,t) \in \{x\} \times (u_1^l, u_{k-1}^r),\\
\phi_{T_2} & (x,t) \in \{x\} \times (u_k^l, u_{k}^r),\\
\{-\alpha/(u_k^l - u_{k-1}^r),\frac{(\tilde\phi^x)^2}{4} - \beta(t-g)^2\}  & (x,t) \in \{x\} \times (u_{k-1}^r, u_{k}^l),\\
0 & \mbox{ otherwise}.
\end{array}
\right.
\end{displaymath}
Let us prove that $\tilde \phi \in K_{\{x\}}$. \\
\begin{eqnarray*}
\left|\int_{t_1}^{t_2} \tilde \phi(x,t) \, dt\right| &=&  \left|\int_{t_1}^{u^l_{k-1}}\phi^x_{T_1}(x,t)\, dt - \alpha + \int_{u^r_k}^{t_2}\phi^x_{T_2}(x,t)\, dt\right| \\
&=& \left|\alpha\frac{(u_{k-1}^l - t_1)}{(u_{k-1}^l - u_1^l)}  - \alpha  + \alpha\frac{(t_2 - u_{k}^r )}{(u_{k}^r - u_k^l)}\right| \leq \alpha
\end{eqnarray*}  
for every $t_1 \leq u_1^l$, $t_2 \geq u_k^r$. As in all the other cases the computation is similar, then $\tilde \phi \in K_{\{x\}}$. Therefore
\begin{displaymath}
\mathcal{G}(T_1,\{x\}) + \mathcal{G}(T_2,\{x\}) = \int_{\mathcal{M} \cap (\{x\} \times \R)} \langle \tilde \phi, \nu_T\rangle \theta \, d\mathcal{H}^1 \leq \mathcal{G}(T,\{x\}).
\end{displaymath}
On the other hand by convexity
\begin{displaymath}
\mathcal{G}(T,\{x\}) \leq \mathcal{G}(T_1,\{x\}) + \mathcal{G}(T_2,\{x\}) = \int_{\mathcal{M} \cap (\{x\} \times \R)} \langle \tilde \phi, \nu_T\rangle \theta \, d\mathcal{H}^1.
\end{displaymath}
So the thesis follows.
\end{dimo}
We can prove the analogue:

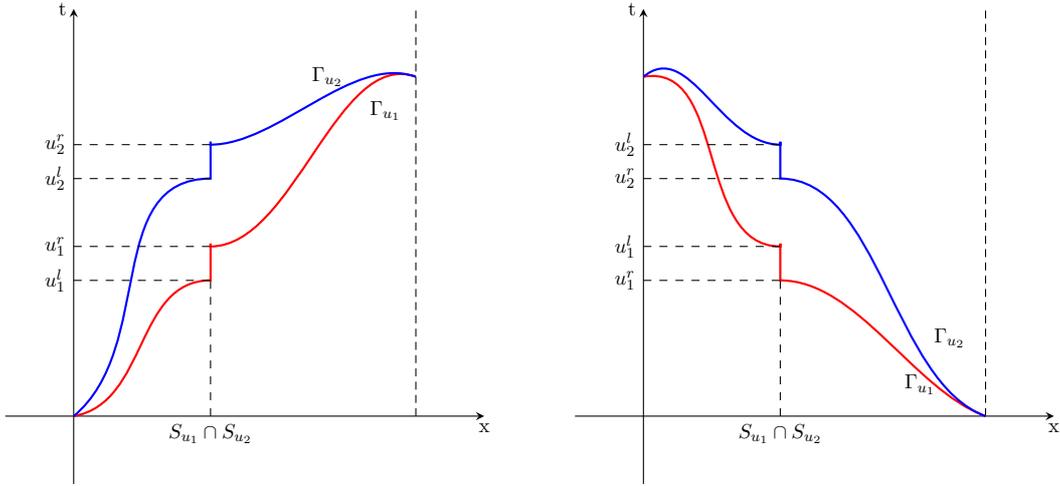
\begin{figure}
\centering
\subfigure{
\begin{tikzpicture}[>=stealth, scale=0.9]
\draw [->](-1,0) -- (6,0) node [scale=0.7][below]{x};
\draw [->](0,-1) -- (0,6) node [scale=0.7][left]{t};
\draw[densely dashed] (5,0) -- (5,6);
\draw  (0,0) [thick][red]
to[out=10,in=180, looseness=1] (2,2) 
to[out=90,in=90, looseness=1] (2,2.5) 
to[out=0,in=160, looseness=0.8] (5,5) node[black, scale=0.7] [below=18pt, left=4pt] {$\Gamma_{u_1}$} ;

\draw  (0,0)[thick][blue]
to[out=40,in=180, looseness=1] (2,3.5)
to[out=90,in=90, looseness=1] (2,4)
to[out=0,in=160, looseness=0.8] (5,5)node[black, scale=0.7] [left=35pt]{$\Gamma_{u_2}$};

\draw  (2,0) node[black, scale=0.7] [below=0.5pt]{$S_{u_1}\cap S_{u_2}$} [dashed] -- (2,2) ;
\draw  (0,2) node[black, scale=0.7] [left=0.5pt]{$u_1^l$} [dashed] -- (2,2) ;
\draw  (0,2.5) node[black, scale=0.7] [left=0.5pt]{$u_1^r$} [dashed] -- (2,2.5) ;
\draw  (0,3.5) node[black, scale=0.7] [left=0.5pt]{$u_2^l$} [dashed] -- (2,3.5) ;

\draw  (0,4) node[black, scale=0.7] [left=0.5pt]{$u_2^r$} [dashed] -- (2,4) ;

\end{tikzpicture}
}\qquad
\subfigure{
\begin{tikzpicture}[>=stealth, scale=0.9]
\draw [->](-1,0) -- (6,0) node [scale=0.7][below]{x};
\draw [->](0,-1) -- (0,6) node [scale=0.7][left]{t};
\draw[densely dashed] (5,0) -- (5,6);
\draw  (0,5) [thick][red]
to[out=10,in=180, looseness=1] (2,2.5) 
to[out=90,in=90, looseness=1]  (2,2)
to[out=0,in=160, looseness=0.8] (5,0) node[black, scale=0.7] [above=42pt, left=7pt] {$\Gamma_{u_2}$} ;

\draw  (0,5)[thick][blue]
to[out=40,in=180, looseness=1]  (2,4)
to[out=90,in=90, looseness=1] (2,3.5)
to[out=0,in=160, looseness=0.8] (5,0)node[black, scale=0.7] [left=35pt, above=8pt]{$\Gamma_{u_1}$};

\draw  (2,0) node[black, scale=0.7] [below=0.5pt]{$S_{u_1}\cap S_{u_2}$} [dashed] -- (2,2) ;
\draw  (0,2) node[black, scale=0.7] [left=0.5pt]{$u_1^r$} [dashed] -- (2,2) ;
\draw  (0,2.5) node[black, scale=0.7] [left=0.5pt]{$u_1^l$} [dashed] -- (2,2.5) ;
\draw  (0,3.5) node[black, scale=0.7] [left=0.5pt]{$u_2^r$} [dashed] -- (2,3.5) ;

\draw  (0,4) node[black, scale=0.7] [left=0.5pt]{$u_2^l$} [dashed] -- (2,4) ;

\end{tikzpicture}
}
\caption{Configuration in Lemma \ref{fillingone} and in Lemma \ref{fillingtwo}}
\label{splitfig}
\end{figure}

\begin{lemma}\label{fillingtwo}
Given $T= \sum_{i=1}^k \Gamma_{u_i}\in C$ such that $u_i$ are ordered in an increasing way. Fix $x\in S_T$ and suppose that we have $u_i^r(x) \leq u_i^l(x)$ for every $i=1 \ldots k$. Suppose in addition that
\begin{displaymath}
u_i^r(x) > u_j^l(x)\quad \mbox{ for every }  \quad i>j.
\end{displaymath}
Then
\begin{displaymath}
\mathcal{G}(T,\{x\}) = \sum_{i=1}^k\mathcal{G}(\Gamma_{u_i},\{x\})= \alpha |\{i: x\in S_{u_i}\}|.
\end{displaymath}
In addition the maximum is achieved and letting $\phi_{T}$ be the vector field realizing the maximum for $T$
\begin{displaymath}
\phi^x_T(x,t) = \alpha/(u_i^r - u_i^l) \quad \mbox{ for every } \quad t\in (u_i^l,u_i^r)
\end{displaymath}
for every $i=1\ldots k$ such that $x\in S_{u_i}$.
\end{lemma}
\begin{dimo}
See Lemma \ref{fillingone}.
\end{dimo}
We are now in position to prove two general statements that are generalizations of Lemmas \ref{fillingone} and \ref{fillingtwo}. 
\begin{prop}\label{summ1}
Consider $T \in C$ such that $T = \sum_{i=1}^k \Gamma_{u_i}$. Fix $x\in S_T$ and suppose that we have $u_i^l(x) \leq u_i^r(x)$ for every $i=1 \ldots k$. Moreover assume that $u^r_{i} (x)  \neq  u^l_{j}(x)$ for every $i,j$ such that $x\in S_{u_i}\cap S_{u_j}$.\\
Then
\begin{displaymath}
\mathcal{G}(T, \{x\}) = \sum_{i=1}^k \mathcal{G}(\Gamma_{u_i}, \{x\}).
\end{displaymath} 
\end{prop}
\begin{dimo}
We can assume without loss of generality that $x\in S_{u_i}$ for every $i=1 \ldots k$.\\
Then it is easy to see that 
\begin{equation}\label{ah}
T\res(\{x\} \times \R) = \sum_{i=1}^{k'} \mu_i[\![\{x\} \times(a_i,a_{i+1})]\!]
\end{equation}
for some $\lambda_i \in \N$ and $a_i \in \R$ with $a_i < a_{i+1}$ and $\lambda_i \neq \lambda_{i+1}$ for every $i$. Let us denote by $\{\lambda_{M_j}\}$ the local maxima of the sequence $\{\lambda_i\}$ and let $\lambda_{m_j}$ be the minimum multiplicity in $\{\lambda_{M_j}, \lambda_{M_j + 1}, \ldots , \lambda_{M_{j+1} - 1}, \lambda_{M_{j+1}}\}$ for every $j$. \\
Thanks to the assumptions, we have that
\begin{equation}\label{maxmin}
k = \sum_j \lambda_{M_j} - \sum_j \lambda_{m_j}.
\end{equation}
Equation \eqref{maxmin} can be proved by induction. Consider $T=\sum_{i=1}^k \Gamma_{u_i}$ associated to a sequence of natural numbers $\{\lambda_i\}_{i=1,\ldots,k'}$ and intervals $(a_i,a_{i+1})$ according to \eqref{ah} with $\lambda_i \neq \lambda_{i+1}$ and $a_i < a_{i+1}$ for every $i$. Let $\Gamma_{w}$ be the graph composing $T$ such that
\begin{displaymath}
\Gamma_w \res (\{x\} \times \R) = [\![\{x\} \times (b,c)]\!]
\end{displaymath}
with $b \neq c$ and $c = a_{k'}$ and call $\hat T = T - \Gamma_w \in C$. We want to show that adding $\Gamma_w$ to $\hat T$ we are increasing the quantity $\sum_j \lambda_{M_j} - \sum_j \lambda_{m_j}$ by one. This fact can be verified considering separately the cases in which $ b\in (a_h,a_{h+1})$ where $\lambda_h$ is a local maximum, a local minimum and none of the two for the sequence $\{\lambda_i\}_{i=1,\ldots,k'}$. 

Then the proof proceeds similarly to the proof of Lemma \ref{fillingone}. 
One can define the following vector field on $\{x\} \times \R$:
\begin{displaymath}
\tilde \phi = \left\{\begin{array}{ll}
\Big(\alpha/(a_{M_{j} + 1}-a_{M_j}),\frac{(\tilde\phi^x)^2}{4} - \beta(t-g)^2\Big) & (x,t) \in \{x\} \times \bigcup_j (a_{M_j}, a_{M_j + 1}),\\
\Big(-\alpha/(a_{m_{j} + 1}-a_{m_{j}}), \frac{(\tilde\phi^x)^2}{4} - \beta(t-g)^2\Big) & (x,t) \in \{x\}\times \bigcup_j  (a_{m_{j}}, a_{m_{j}+1}), \\
0 & \mbox{ otherwise,}
\end{array}
\right.
\end{displaymath}
proving that $\tilde \phi \in K_{\{x\}}$, similarly as in Lemma \ref{fillingone}. Then, thanks to \eqref{maxmin}, one obtains
\begin{displaymath}
\mathcal{G}(T,\{x\}) = \alpha k = \sum_{i=1}^k \mathcal{G}(\Gamma_{u_i}, \{x\})
\end{displaymath}
as we wanted to prove.

\end{dimo}

\begin{prop}\label{summ2}
Consider $T \in C$ such that $T = \sum_{i=1}^k \Gamma_{u_i}$. Fix $x\in S_T$ and suppose that we have $u_i^r(x) \leq u_i^l(x)$ for every $i=1 \ldots k$. Moreover assume that $u^l_{i} (x)  \neq  u^r_{j}(x)$ for every $i,j$ such that $x\in S_{u_i}\cap S_{u_j}$.\\
Then
\begin{displaymath}
\mathcal{G}(T, \{x\}) = \sum_{i=1}^k \mathcal{G}(\Gamma_{u_i}, \{x\}).
\end{displaymath} 
\end{prop}
\begin{dimo}
See Proposition \ref{summ1}.
\end{dimo}
Now Theorem \ref{split} is an immediate consequence of the previous propositions.
\vspace{2mm}\\
\emph{Proof of Theorem \ref{split}}
\vspace{2mm}\\
Fix $x\in S_T$ and define
\begin{displaymath}
\mathcal{I} = \{i=1\ldots k: u^l_{i} (x) \leq u^r_{i}(x)\} \qquad \mathcal{J} = \{i=1\ldots k: u^l_{i} (x) > u^r_{i}(x)\}
\end{displaymath}
and call $T_\mathcal{I} = \sum_{i\in \mathcal{I}}\Gamma_{u_i}$ and $T_\mathcal{J} = \sum_{i\in \mathcal{J}}\Gamma_{u_i}$. 
Moreover let $\phi_{\mathcal{I}}$ ($\phi_\mathcal{J}$) be the vector field realizing the maximum in $\mathcal{G}(T_\mathcal{I},\{x\})$ ($\mathcal{G}(T_\mathcal{J}),\{x\})$. From Proposition \ref{summ1} and \ref{summ2} it is easy to see that $\phi^x_\mathcal{\mathcal{I}} \leq 0$ outside the support of $T_\mathcal{I}$ restricted to $\{x\} \times \R$ and $\phi^x_\mathcal{J} \geq 0$ outside the support of $T_\mathcal{I}$ restricted to $\{x\} \times \R$. Therefore defining $\tilde \phi =  \phi_\mathcal{I} + \phi_\mathcal{J}$, as we assumed that there is no cancellation on the jumps by Lemma \ref{nocanc}, we have that $\tilde \phi \in K_{\{x\}}$ and
\begin{displaymath}
\mathcal{G}(T_\mathcal{I},\{x\}) + \mathcal{G}(T_\mathcal{J},\{x\}) = \int_{\{x\} \times \R}  \langle\phi_\mathcal{I}^x + \phi_\mathcal{J}^x,\nu_T\rangle \, d\|T\|  \leq \mathcal{G}(T,\{x\}).
\end{displaymath} 
So by convexity 
\begin{displaymath}
\mathcal{G}(T_\mathcal{I},\{x\}) + \mathcal{G}(T_\mathcal{J},\{x\}) = \mathcal{G}(T,\{x\}).
\end{displaymath}
Finally we apply Proposition \ref{summ1} and \ref{summ2} to $T_\mathcal{I}$ and $T_\mathcal{J}$ to get the thesis.
\begin{flushright}
$\Box$
\end{flushright} 
We conclude this section with a lemma that shows that we can reduce any combination of graphs belonging to $C$ to a combination of graphs, all with the same multiplicity. We are going to use this property in the proof of the convex decomposition formula in the next section.
\begin{lemma}\label{samemult}
Consider $T_1,T_2\in C$ and $x \in S_{T_1} \cap S_{T_2}$. Suppose that $T_1\res (\{x\} \times \R) = \sum_{i=1}^k \lambda_i [\![\{x\}\times (a_i,a_{i+1})]\!]$ with $a_i < a_{i+1}$ and let $\{M_j\}_{j\in J}$ be the indexes of the maximums of the multiplicities. Assume in addition that $T_2 \res (\{x\} \times \R) = \nu \sum_{j\in J} [\![\{x\}\times (a_{M_j}, a_{M_j+1})]\!]$ for some $\nu>0$. Then we have
\begin{equation}\label{kin}
\mathcal{G}(T_1 + T_2,\{x\}) = \mathcal{G}(T_1,\{x\}) + \mathcal{G}(T_2,\{x\}).
\end{equation}  
\end{lemma}
\begin{dimo}
Consider the vector field $\tilde \phi \in K_{\{x\}}$ defined in Lemma \ref{summ1}: 
\begin{displaymath}
\tilde \phi = \left\{\begin{array}{ll}
\Big(\alpha/(a_{M_{j} + 1}-a_{M_j}),\frac{(\tilde\phi^x)^2}{4} - \beta(t-g)^2\Big) & (x,t) \in \{x\} \times \bigcup_j (a_{M_j}, a_{M_j + 1}),\\
\Big(-\alpha/(a_{m_{j} + 1}-a_{m_{j}}), \frac{(\tilde\phi^x)^2}{4} - \beta(t-g)^2\Big) & (x,t) \in \{x\}\times \bigcup_j  (a_{m_{j}}, a_{m_{j}+1}), \\
0 & \mbox{ otherwise}.
\end{array}
\right.
\end{displaymath}
Thanks to \eqref{maxmin} we have
\begin{displaymath}
\mathcal{G}(T_1, \{x\}) = \int_{\mathcal{M}_1 \cap (\{x\} \times \R)} \langle \tilde \phi^x, \nu_{T_1}\rangle \theta_{T_1}\, d\Ha^1 \quad \mbox{and} \quad \mathcal{G}(T_2, \{x\}) = \int_{\mathcal{M}_2 \cap (\{x\} \times \R)} \langle \tilde \phi^x, \nu_{T_2}\rangle\theta_{T_2}\, d\Ha^1,
\end{displaymath}
where $T_i = (\mathcal{M}_i, \theta_{T_i}, \nu_{T_i})$ for $i=1,2$. Hence setting $T_1 + T_2 =(\mathcal{M},\theta,\nu)$ we have
\begin{eqnarray*}
\mathcal{G}(T_1, \{x\}) + \mathcal{G}(T_2, \{x\}) & = & \int_{\mathcal{M}_1 \cap (\{x\} \times \R)} \langle \tilde{\phi}^x, \nu_{T_1}\rangle \theta_{T_1}\, d\Ha^1 + \int_{\mathcal{M}_2 \cap (\{x\} \times \R)} \langle \tilde \phi^x, \nu_{T_2}\rangle\theta_{T_2}\, d\Ha^1\\
& = & \int_{\mathcal{M} \cap (\{x\} \times \R)} \langle \tilde \phi^x, \nu \rangle \theta \, d\Ha^1 \leq G(T_1 + T_2).
\end{eqnarray*}
As the opposite inequality follows by convexity, we infer \eqref{kin}.

\end{dimo}

\begin{cor}\label{same}
Given $T_1 = \sum_{i=1}^{k} \lambda_i \Gamma_{u_i}$, let $\{M_j\}_{j\in J}$ be the indexes of the maximums of the multiplicities. Given $T_2 =  \nu\sum_{i \in J} \Gamma_{u_i}$ with $\nu > 0$ we have that $G(T_1 + T_2) = G(T_1) + G(T_2)$. 
\end{cor}
\begin{dimo}
Notice that by Lemma \ref{overlapping} it is enough to prove the thesis for every $x\in S_{T_2}\cap S_{T_1}$. Thanks to Lemma \ref{samemult} one has
\begin{displaymath}
G(T_1 + T_2, \{x\}) =  G(T_1, \{x\}) + G(T_2, \{x\}).
\end{displaymath}
\end{dimo}
\subsection{Convex decomposition formula}
As anticipated in the introduction, this section is devoted to the proof of a decomposition formula for the Mumford-Shah functional in one dimension. This formula resembles closely a generalized coarea formula for functionals and it is performed for a finite combination of graphs with multiplicity. It is interesting to notice that the counterexample in the end of Remark \ref{coareafalse} is ``solved'' by this decomposition, but it is difficult to generalize it to the continuous case. However it gives a strong indication on how this decomposition should be performed at least in dimension one. The higher dimensional case is a completely different issue, as the convex decomposition formula we are going to present strongly relies on the one dimensional structure of the problem and cannot be extended in an easy way.
\begin{figure}
\centering
\subfigure{
\begin{tikzpicture}[>=stealth, scale=0.9]
\draw [->](-1,0) -- (6,0) node [scale=0.7][below]{x};
\draw [->](0,-1) -- (0,6) node [scale=0.7][left]{t};
\draw[densely dashed] (5,0) -- (5,6);
\draw  (0,0) [thick][red]
to[out=10,in=180, looseness=1] (2,2) 
to[out=90,in=90, looseness=1] (2,3) 
to[out=0,in=160, looseness=0.8] (5,5) node[black, scale=0.7] [below=18pt, left=4pt] {$\Gamma_{u_1}$} ;

\draw  (0,0)[thick][blue]
to[out=40,in=180, looseness=1] (2,3)
to[out=90,in=-90, looseness=1] (2,4)
to[out=0,in=160, looseness=0.8] (5,5)node[black, scale=0.7] [left=35pt]{$\Gamma_{u_2}$};

\draw  (2,0) node[black, scale=0.7] [below=0.5pt]{$S_{u_1}\cap S_{u_2}$} [dashed] -- (2,2) ;
\draw  (0,2) node[black, scale=0.7] [left=0.5pt]{$u_1^l$} [dashed] -- (2,2) ;
\draw  (0,3) node[black, scale=0.7] [left=0.5pt]{$u_1^r = u_2^l$} [dashed] -- (2,3) ;

\draw  (0,4) node[black, scale=0.7] [left=0.5pt]{$u_2^r$} [dashed] -- (2,4) ;

\end{tikzpicture}
} \qquad
\subfigure{
\begin{tikzpicture}[>=stealth, scale=0.9]
\draw [->](-1,0) -- (6,0) node [scale=0.7][below]{x};
\draw [->](0,-1) -- (0,6) node [scale=0.7][left]{t};
\draw[densely dashed] (5,0) -- (5,6);
\draw  (0,0) [thick][green]
to[out=10,in=180, looseness=1] (2,2) 
to[out=90,in=-90, looseness=1] (2,4) 
to[out=0,in=160, looseness=1] (5,5) node[black, scale=0.7] [below=18pt, left=4pt] {$\Gamma_{w_2}$} ;

\draw  (0,0)[thick][violet]
to[out=40,in=180, looseness=1] (2,3)
to[out=0,in=160, looseness=0.8] (5,5)node[black, scale=0.7] [left=35pt]{$\Gamma_{w_1}$};

\draw  (2,0) node[black, scale=0.7] [below=0.5pt]{$S_{w_1}$} [dashed] -- (2,2) ;
\draw  (0,2) node[black, scale=0.7] [left=0.5pt]{$w_1^l$} [dashed] -- (2,2) ;

\draw  (0,4) node[black, scale=0.7] [left=0.5pt]{$w_1^r$} [dashed] -- (2,4) ;

\end{tikzpicture}

}

\caption{Convex decomposition of two SBV graphs}
\label{space}
\end{figure}
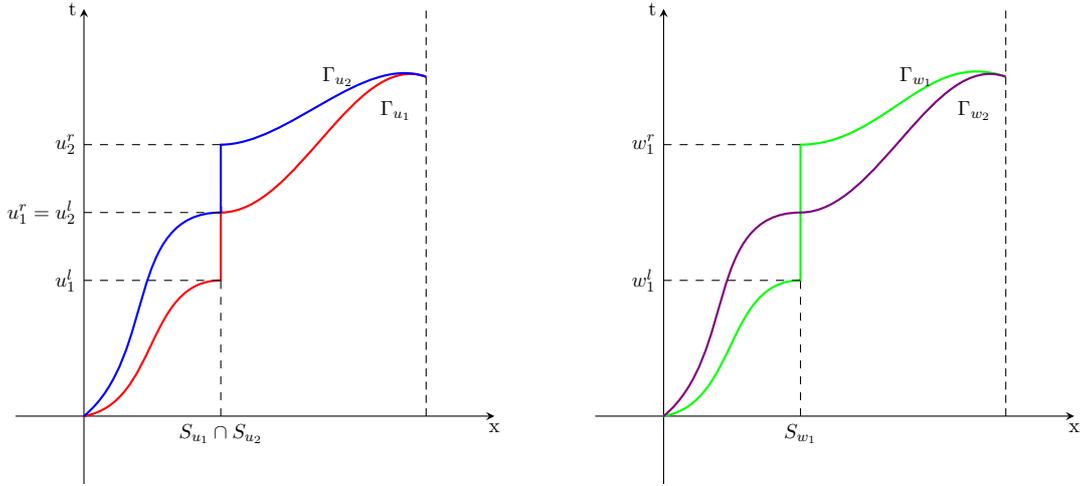

\begin{prop}\label{coareasame}
Given $T= \sum_{i=1}^k \Gamma_{u_i} \in C$ such that $|S_T| <+\infty$ there exists $\{w_i\}_{i=1\ldots k} \subset SBV(I)$ such that $T = \sum_{i=1}^k\Gamma_{w_i}$ and 
\begin{displaymath}
G(T) =  \sum_{i=1}^k G(\Gamma_{w_i}).
\end{displaymath}
\end{prop}
\begin{dimo}
Set $I=(a,b)$ and $S_T = \{p_1,\ldots,p_N\} \subset I$ (the cardinality of $S_T$ is finite by assumption).
Choose $N-1$ points $\{x_s\}$ so that $p_s < x_s < p_{s+1}$ and define a partition of $(a,b]$ as 
\begin{displaymath}
(a,b]=\bigcup_{s=0}^{N-1} (x_s,x_{s+1}]
\end{displaymath}
where $x_0 = a$ and $x_N =b$.\\
Consider an interval $(x_s,x_{s+1})$ and define $\mathcal{I} \subset \{1,\ldots,k\}\times \{1,\ldots,k\}$ as
\begin{equation}
\mathcal{I} := \{(i,j) : p_s \in S_{u_i} \cap S_{u_j} \ \mbox{and} \  u^r_i(p_s) = u^l_j(p_s) \}.
\end{equation}

For a given $(i,j) \in \mathcal{I}$ define the following functions (see Figure \ref{space}):

\begin{displaymath}
w_i = \left\{
\begin{array}{ll}
u_{i} & \mbox{ for } x \in (x_s,p_s)\\
u_{j} & \mbox{ for } x \in (p_s,x_{s+1})
\end{array}
\right.
\end{displaymath}
and
\begin{displaymath}
w_j = \left\{
\begin{array}{ll}
u_{j} & \mbox{ for } x \in (x_s,p_s)\\
u_{i} & \mbox{ for } x \in (p_s,x_{s+1}).
\end{array}
\right. 
\end{displaymath}
We have that $\Gamma_{u_i} + \Gamma_{u_j} = \Gamma_{w_i} + \Gamma_{w_j}$ in $(x_s,x_{s+1})$ and $w_i^r(p_s) \neq w_j^l(p_s)$. Then, considering the new collection of functions with $u_i,u_j$ substituted with $w_i$ and $w_j$, we can repeat this operation. It is easy to see that we can perform this procedure only a finite number of times (until the set $\mathcal{I}$ is empty), as at every step we are strictly decreasing the cardinality of $\mathcal{I}$ by at least one. In this way we produce a family of functions $\{w_i\} \subset SBV((x_s,x_{s+1}))$ such that
$\sum_{i=1}^k \Gamma_{w_i} = \sum_{i=1}^k \Gamma_{u_i}$ in $(x_s, x_{s+1})$ and $w^r_{i} (p_s) \neq w^l_{j}(p_s)$ for every $(i,j)$ such that $p_s \in S_{w_i}\cap S_{w_j}$. \\
Finally we can construct a family of functions $\{w_i\} \subset SBV(I)$ repeating this procedure inductively starting from the first interval $(a,x_1)$ and ending in $(x_{n-1},b)$. The family $\{w_i\} \subset SBV(I)$ has the following properties: 
\begin{displaymath}
\sum_{i=1}^k \Gamma_{w_i} = \sum_{i=1}^k \Gamma_{u_i} \quad \mbox{and} \quad w^r_{i} (x) \neq w^l_{j}(x)\  \  \forall x,i,j \mbox{ such that } x \in S_{w_i}\cap S_{w_j}.
\end{displaymath}
Hence using Theorem \ref{split} and Proposition \ref{decomposition} one obtains the thesis.

\end{dimo}

\begin{thm}[Convex decomposition formula]\label{coarea}
Given $T= \sum_{i=1}^k \lambda_i\Gamma_{u_i}$ such that $|S_T| < +\infty$ there exists $k' \in \mathbb{N}$, $\{\mu_i\}_{i=1\ldots k'} \geq 0$ and $\{w_i\}_{i=1\ldots k'} \subset SBV(I)$ such that $T = \sum_{i=1}^{k'} \mu_i\Gamma_{w_i}$ and 
\begin{equation}\label{coar}
G(T) =  \sum_{i=1}^{k'} \mu_i G(\Gamma_{w_i}).
\end{equation}
\end{thm}
\begin{dimo}
Consider $T= \sum_{i=1}^k \lambda_i \Gamma_{u_i} \in C$ and suppose without loss of generality that also $\lambda_i$ are ordered and $\lambda_k$ is the maximum (if the multiplicities are not ordered the proof is analogous). Then $T$ can be rewritten as
\begin{displaymath}
T= (\lambda_k-\lambda_{k-1}) \Gamma_{u_k} + \lambda_{k-1} \Gamma_{u_k} + \sum_{i=1}^{k-1} \lambda_i \Gamma_{u_i}.
\end{displaymath}
Hence by Corollary \ref{same}
\begin{displaymath}
G(T) = G((\lambda_k-\lambda_{k-1}) \Gamma_{u_k}) + G\left(\lambda_{k-1} \Gamma_{u_k} + \sum_{i=1}^{k-1} \lambda_i \Gamma_{u_i}\right).
\end{displaymath}
Then one can rewrite
\begin{displaymath}
\lambda_{k-1} \Gamma_{u_k} + \sum_{i=1}^{k-1} \lambda_i \Gamma_{u_i} = \lambda_{k-2} (\Gamma_{u_k} + \Gamma_{u_{k-1}}) +  (\lambda_{k-1} -\lambda_{k-2}) (\Gamma_{u_k} +  \Gamma_{u_{k-1}}) + \sum_{i=1}^{k-2} \lambda_i \Gamma_{u_i}
\end{displaymath}
and applying again Corollary \ref{same}
\begin{eqnarray}
G\left(\lambda_{k-1} \Gamma_{u_k} + \sum_{i=1}^{k-1} \lambda_i \Gamma_{u_i}\right) &=& 
G((\lambda_{k-1} -\lambda_{k-2}) (\Gamma_{u_k} +  \Gamma_{u_{k-1}}) ) \nonumber\\
&+& G\left(\lambda_{k-2} (\Gamma_{u_k} + \Gamma_{u_{k-1}}) + \sum_{i=1}^{k-2} \lambda_i \Gamma_{u_i}\right).\label{tr}
\end{eqnarray}
By Proposition \ref{coareasame} there exists $u^2_k$ and $u^2_{k-1}$ SBV functions such that $\Gamma_{u^2_k} +  \Gamma_{u^2_{k-1}} = \Gamma_{u_k} + \Gamma_{u_{k-1}}$ and
\begin{displaymath}
\eqref{tr}= G((\lambda_{k-1} -\lambda_{k-2}) \Gamma_{u^2_k}) + G((\lambda_{k-1} -\lambda_{k-2}) \Gamma_{u^2_{k-1}})  +  G\left(\lambda_{k-2} (\Gamma_{u_k} + \Gamma_{u_{k-1}}) + \sum_{i=1}^{k-2} \lambda_i \Gamma_{u_i}\right)
\end{displaymath}
and so on. Repeating this procedure $k$ times one gets to
\begin{displaymath}
G(T) = \sum_{i=2}^{k}\sum_{j=i}^k (\lambda_{i} -\lambda_{i-1})G(\Gamma_{u^{k-i+1}_j}) + G\left(\sum_{i=1}^k \lambda_1 \Gamma_{u_i}\right).
\end{displaymath}
Hence, applying again Proposition \ref{coareasame} to the last term we obtain the desired decomposition \eqref{coar}.

\end{dimo}
\section{Existence of calibration as a functional defined on currents}\label{existence}
We now want to show an application of the previous convex decomposition formula to the existence of calibration for the Mumford-Shah type functionals. Firstly we set the Dirichlet problem in $I'$ associated to the previous functional $G$.
We recall that $G$ takes values in the convex cone defined as
\begin{equation}
C:= \left\{T= \sum_{i=1}^{k} \lambda_i \Gamma_{u_i} : k\in \mathbb{N}, \lambda_i \in \mathbb{R}_+, u_i \in SBV(I)\right\}.
\end{equation}
Consider $S\in C$ and setting $I_{\R} := (I\setminus I') \times \R$ define
\begin{displaymath}
\psi_{G}(S) = \inf\{G(T) : T\in C, \ T\res I_{\R} = S \res I_{\R}\}.
\end{displaymath}

\begin{prop}\label{conv}
The functional $\psi_{G}$ is convex in $C$.
\end{prop}
\begin{dimo}
As $G$ is convex and the constraint is linear the proof is straightforward.
\end{dimo}
It is easy to see that by the convex decomposition formula in Theorem \ref{coarea} we have the following theorem:
\begin{thm}\label{mainthm1}
If $u \in SBV(I)$ is a Dirichlet minimizer in $I'$ of $F$, then $\psi_{G}(\Gamma_u) = G(\Gamma_u) = F(u)$.
\end{thm}
\begin{dimo}
Consider $T = \sum_{i=1}^k \lambda_i \Gamma_{u_i} \in C$ such that $T\res I_{\R} = \Gamma_u \res I_{\R}$. 
Without loss of generality we can suppose that $|S_T| <+\infty$. 
%Then letting $\pi:\R^2 \rightarrow \R$ the projection on the first component we have
%\begin{equation}\label{push}
%\partial [\![I']\!] = \pi^\# (\partial \Gamma_u) = \pi^\# (\partial T) = (\partial [\![I']\!]) \sum_{i=1}^k  \lambda_i.  
%\end{equation}
%Hence $\sum_{i=1}^k \lambda_i = 1$. 
Then by Theorem \ref{coarea} there exist $k'$ and $\{\mu_i\}_{i=1,\ldots, k'}>0$ such that
\begin{equation}\label{lastcor}
G(T) = G\left(\sum_{i=1}^k \lambda_i \Gamma_{u_i}\right) = \sum_{i=1}^{k'}\mu_i G(\Gamma_{w_i}) = \sum_{i=1}^{k'}\mu_i F(w_i)
\end{equation}
and $\sum_{i=1}^k \lambda_i \Gamma_{u_i} = \sum_{i=1}^{k'}\mu_i \Gamma_{w_i}$. As $T\res I_{\R} = \Gamma_u \res I_{\R}$, we infer that $\sum_{i=1}^{k'}\mu_i = 1$ and
$w_i = u$ in $I\setminus I'$ for every $i=1,\ldots,k'$. Finally, from Formula \eqref{lastcor}, using that $u$ is a Dirichlet minimizer of $F$ in $I'$ we obtain that $\psi_{G}(\Gamma_u) = G(\Gamma_u)$.

\end{dimo}

\begin{rmk}
Theorem \ref{mainthm1} can be obtained also with similar techniques as the ones presented in \cite{CRF}, provided that the functional $F$ satisfies stronger regularity assumptions (see Remark \ref{gener}).
\end{rmk}

The previous theorem allows us to state a weak existence results for calibrations as an application of Hahn-Banach theorem. Let 
\begin{displaymath}
\hat C = \left\{T = \sum_{i=1}^k \lambda_i \Gamma_{u_i}: k\in \N, \lambda_i \in \R, u_i \in SBV(I)\right\}
\end{displaymath}
be the double cone and denote by $Hom(\hat C)$ the set of all the linear maps from $\hat C$ to $\R$. 
We define the following notion of calibration for linear combinations of graphs:
\begin{defi}[Calibration for minimal graphs]\label{calcur}
Given $u \in SBV(I)$ and $\Gamma_u$ its associated graph, we say that $\xi\in \mbox{Hom}(\hat C)$ is a calibration for $\Gamma_u$ with respect to $G$ if
\begin{itemize}
\item[i)] $\xi(\Gamma_u) = G(\Gamma_u) = F(u)$,
\item[ii)] $\xi(T) = 0$ for every $T \in \hat C$ such that $T\res I_\R = 0$,
\item[iii)] $\xi (T) \leq \hat G(T)$ for every $T \in \hat C$,
\end{itemize}
where $\hat G : \hat C \rightarrow \R$ is the extension of $G$ to $\hat C$ according to Formula \eqref{operator}. 
\end{defi}
\begin{thm}
Given $u\in SBV(I)$ a Dirichlet minimizer of $F$ in $I'$ there exists a calibration for $\Gamma_u$ with respect to $G$ according to Definition \ref{calcur}.
\end{thm}
\begin{dimo}
From Theorem \ref{mainthm1} follows that
\begin{displaymath}
G(\Gamma_{u}) = \psi_{G}(\Gamma_u). 
\end{displaymath}
We firstly notice that, as a consequence of the definition of $K$, we have that
\begin{equation*}
\hat G(T) = \left\{
\begin{array}{ll}
G(T) & \mbox{if } T\in C\\
+\infty & \mbox{of }T\in \hat C \setminus C.
\end{array}
\right.
\end{equation*}
We define $\psi_{\hat G} : \hat C \rightarrow \R$ as  
\begin{equation*}
\psi_{\hat G}(S) := \inf\{\hat G(T) : T\in \hat C, \ T\res I_{\R} = S \res I_{\R}\},
\end{equation*}
that is convex and such that $\psi_{\hat G}(\Gamma_{u}) = G(\Gamma_{u}) > 0$.\\
Consider the vector subspace $L = \{a\Gamma_u: a\in \R\}$ and define $\psi: L \rightarrow \R$ as $\psi(a\Gamma_u) = a\psi_{\hat G}(\Gamma_u)$ clearly linear. As we have that $\psi \leq \psi_{\hat G}$ on $L$, by Hahn-Banach theorem there exists $\xi\in Hom(\hat C, \R)$ such that 
\begin{equation}\label{too}
\xi(\Gamma_{u}) =  \psi(\Gamma_u) = \psi_{\hat G}(\Gamma_{u}) \quad \mbox{ and }\quad  \xi(T) \leq \psi_{\hat G}(T) \ \quad \forall T\in \hat C.
\end{equation}
We want to prove that $\xi$ is a calibration according to Definition \ref{calcur}. Let $T_0\in \hat C$ be such that $T_0\res I_\R  = 0$, then 
\begin{displaymath}
\psi_{\hat G}(T_0)= \inf\{\hat G(S) : S\in \hat C, \  S\res I_\R = 0\} \leq G(0) = 0.
\end{displaymath}
In combination with \eqref{too} this implies $\xi(T)\leq 0$ for every $T_0 \in \hat C$ such that $T_0\res I_\R  = 0$.\\
So, as $\xi$ is an homeomorphism, one has also that $\xi(T_0)= 0$, so that $(ii)$ holds. 
Moreover from \eqref{too}, $\xi(\Gamma_u) = \psi_{\hat G}(\Gamma_{u}) = F(u)$ that is $(i)$.\\
Let us show that also $(iii)$ is satisfied: if $T\in \hat C \setminus C$ then $G(T) = +\infty$ and so there is nothing to prove. On the other hand given $T = \sum_{i=1}^{k} \lambda_i \Gamma_{u_i} \in C$ with $\lambda_i \in \R_+$ by \eqref{too} and using the definition of $\psi_{\hat G}$
\begin{displaymath}
\xi(T) \leq \psi_{\hat G}(T) \leq \hat G(T).
\end{displaymath}
Hence $\xi$ is a calibration according to Definition \ref{calcur}.
\end{dimo}

\bibliographystyle{plain}
\bibliography{biblio}

\end{document}